\documentclass[a4paper,12pt,reqno]{amsart}

\usepackage{rotating,pdflscape,color,amssymb,subfigure,psfrag,amsmath,eufrak,bbm,epsfig}
\usepackage{a4wide,amscd}
\usepackage{tikz} % doc : chercher PGF
 \usepackage[usenames,dvipsnames]{pstricks}
 \usepackage{pst-grad} % For gradients
 \usepackage{pst-plot} % For axes

\definecolor{vdarkred}{rgb}{0.6,0,0.2}
\definecolor{vdarkblue}{rgb}{0,0.2,0.6}
\usepackage[pdftex, colorlinks, linkcolor=vdarkblue,citecolor=vdarkred]{hyperref}
\usepackage[small]{caption}

\DeclareMathOperator{\im}{Im}

\newcommand{\ld}{\ldots}
\newcommand{\beg}{\begin}
\newcommand{\en}{\end}

\newcommand{\trm}{\textrm}
 
\newcommand{\bgt}{\begin{itemize}}
\newcommand{\ent}{\end{itemize}}
\newcommand{\ite}{\item}

\newcommand{\op}{\operatorname}
\newcommand{\eqre}{\eqref}
\newcommand{\re}{\ref}
\newcommand{\la}{\label}

\newcommand{\si}{\sigma}
\newcommand{\Sy}{\mathfrak{S}}

\newcommand{\diag}{\operatorname{diag}}

 \newcommand{\bgn}{\begin{enumerate}}
\newcommand{\enn}{\end{enumerate}}
\newcommand{\ds}{\displaystyle}
 
\newcommand{\p}{\mathbb{P}}

\newcommand{\supp}{\operatorname{supp}}

\newcommand{\Tr}{\operatorname{Tr}}

\newcommand{\ninf}{\underset{n\to\infty}{\longrightarrow}}

\newcommand{\E}{\mathbb{E}}

\newcommand{\R}{\mathbb{R}}
\newcommand{\C}{\mathbb{C}}

\newcommand{\ud}{\mathrm{d}}

\newcommand{\pro}{probability }

\newcommand{\f}{\frac}
\newcommand{\ff}{\frac{1}}
\newcommand{\lf}{\left}
\newcommand{\ri}{\right}

\newcommand{\st}{such that }

\newcommand{\lam}{\lambda}

\newcommand{\ti}{\times}

\newcommand{\vfi}{\varphi}
\newcommand{\ste}{\, ;\, }

\newcommand{\eps}{\varepsilon}

\newcommand{\A}{\mathbf{A}}

\newcommand{\bck}{\backslash}
\newcommand{\al}{\alpha}

\newcommand{\eqlaw}{\stackrel{\textrm{law}}{=}}

\newcommand{\ovl}{\overline}

\newcommand{\bbm}{\begin{bmatrix}}
\newcommand{\ebm}{\end{bmatrix}}
\newcommand{\bes}{\begin{equation*}}
\newcommand{\ees}{\end{equation*}}
\newcommand{\be}{\begin{equation}}
\newcommand{\ee}{\end{equation}}
\newcommand{\beqy}{\begin{eqnarray}}
\newcommand{\eeqy}{\end{eqnarray}}
\newcommand{\beq}{\begin{eqnarray*}}
\newcommand{\eeq}{\end{eqnarray*}}
\newcommand{\one}{\mathbbm{1}}
\newcommand{\lto}{\longrightarrow}

\newcommand{\ie}{\emph{i.e. }}

\newcommand{\bpm}{\begin{pmatrix}}
\newcommand{\epm}{\end{pmatrix}}

\newcommand{\cd}{\cdots}

\newcommand{\wg}{\op{Wg}}

\newcommand{\bpr}{\beg{proof}}
\newcommand{\epr}{\en{proof}}

\newcommand{\del}{\delta}

\newcommand{\ka}{\kappa}

\newcommand{\bi}{\mathbf{i}}

\newcommand{\bA}{\mathbf{A}}

\newcommand{\bV}{\mathbf{V}}
\newcommand{\bU}{\mathbf{U}}
\newcommand{\bT}{\mathbf{T}}

\newcommand{\bM}{\mathbf{M}}

\newcommand{\bj}{\mathbf{j}}

\newtheorem{Th}{Theorem}[]

\newtheorem{propo}[Th]{Proposition}

\theoremstyle{definition}

 \allowdisplaybreaks
 
\date{\today}
\title[]{Exponential bounds for the support convergence in the Single Ring Theorem}
\author[]{Florent Benaych-Georges} \address{MAP 5, UMR CNRS 8145 - Universit\'e Paris Descartes, 45 rue des Saints-P\`eres 75270 Paris cedex~6,  France.} \email{florent.benaych-georges@parisdescartes.fr}



\begin{document}
\maketitle

\begin{abstract}
We consider an $n\ti n$  matrix of the form $\bA=\bU\bT\bV$, with $\bU, \bV$ some independent Haar-distributed unitary matrices and $\bT$ a deterministic matrix. We prove that for $k\sim n^{1/6}$ and $b^2:=\frac{1}{n}\Tr(|\bT|^2)$, as $n$ tends to infinity, we have  $$\E\Tr (\bA^{k}(\bA^{k})^*)   \ \lesssim \  b^{2k}\qquad\trm{ and }\qquad\E[|\Tr (\bA^{k})|^2]   \lesssim   b^{2k}.$$
This  gives a simple proof (with slightly weakened hypothesis) of the convergence of the support in the Single Ring Theorem, improves the available error bound for this convergence from $n^{-\alpha}$ to $e^{-cn^{1/6}}$ and  proves that  the rate of this convergence is at most $n^{-1/6}\log n$.
\en{abstract}

%arxiv abstract : We consider an $n$ by $n$ matrix of the form $A=UTV$, with $U, V$ some independent Haar-distributed unitary matrices and $T$ a deterministic matrix. We prove that for $k\sim n^{1/6}$ and $b^2:=\frac{1}{n}\operatorname{Tr}(|T|^2)/n$, as $n$ tends to infinity, we have $$\mathbb{E} \operatorname{Tr} (A^{k}(A^{k})^*) \ \lesssim \ b^{2k}\qquad \textrm{ and } \qquad\mathbb{E}[|\operatorname{Tr} (A^{k})|^2] \ \lesssim \ b^{2k}.$$ This gives a simple proof (with slightly weakened hypothesis) of the convergence of the support in the Single Ring Theorem, improves the available error bound for this convergence from $n^{-\alpha}$ to $e^{-cn^{1/6}}$ and proves that the rate of this convergence is at most $n^{-1/6}\log n$.

%In version v3, we added a corollary (Corollary 4) weakening the hypotheses in the convergence of the support in the Single Ring Theorem.

\section{Introduction}
The Single Ring Theorem, by Guionnet, Krishnapur and Zeitouni \cite{GUI}, describes the   empirical distribution of the eigenvalues of a large generic matrix with prescribed singular values, \ie an $n\ti n$  matrix of the form $\bA=\bU\bT\bV$, with $\bU, \bV$ some independent Haar-distributed unitary matrices and $\bT$ a deterministic matrix whose singular values are the ones prescribed.    More precisely, under some technical hypotheses\footnote{These hypotheses have    have been weakened by  Rudelson and Vershynin in \cite{RUD} and by  Basak and  Dembo in \cite{BD13}.}, as the dimension $n$ tends to infinity, if  the empirical distribution of the  singular values of $\bA$ converges to a compactly supported limit measure $\Theta$ on the real line, 
then the empirical eigenvalues distribution of $\bA$ converges to a limit measure $\mu$ on the complex plane
which depends only on $\Theta$. The limit measure $\mu$ (see Figure \re{figure_intro})
\begin{figure}[ht]
\centering
\includegraphics[scale=.6]{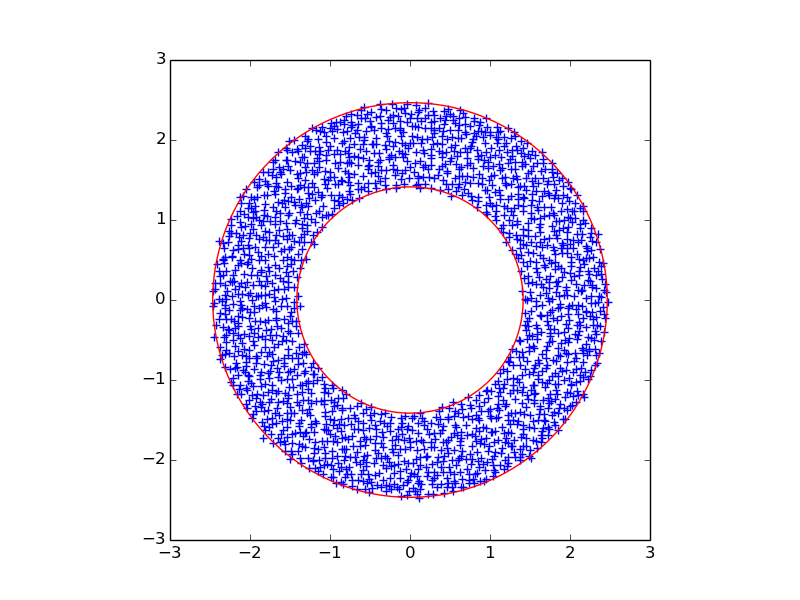}
\caption{Spectrum of $\bA$  when the $s_i$'s are uniformly distributed on $[0.5,4]$, so that $a\approx 1.41$ and $b\approx 2.47$ (here,  $n=2.10^3$).}\label{figure_intro}
\end{figure} 
is rotationally invariant in $\C$ and its support   is the annulus $\{z\in \C\ste a\le |z|\le b\}$, with $a,b\ge 0$ \st \be\la{9914}a^{-2}=\int x^{-2}\ud\Theta( x)\qquad \trm{ and }\qquad b^{2}=\int x^{2}\ud\Theta(  x).\ee In \cite{GUI2}, Guionnet and Zeitouni  also  proved the convergence in \pro of the support of  the empirical eigenvalues distribution of $\bA$ to the support of $\mu$. 
The reason why the radii $a$ and $b$ of the borders of the support of $\mu$ are given by \eqre{9914} is related to the earlier  work \cite{haag2} by Haagerup and Larsen about $R$-diagonal elements in free probability theory but
has no 	simple explanation: the matrix $\bA$ is far from being normal, hence its spectral radius should be  smaller  than its operator norm, \ie than %the supremum of the support of the measure 
  the $L^\infty$-norm\footnote{To be precise, we should say the ``$L^\infty$-norm of \emph{a 
 $\Theta$-distributed r.v.}" rather than ``$L^\infty$-norm of $\Theta$". The same is true for the $L^2$-norm hereafter.} of  $\Theta$, but, up to our knowledge, there is no evidence why this modulus has to be close to the $L^2$-norm of $\Theta$, as follows from \eqre{9914}. 
 
 Another way to see the problem is the following one. In \cite{RUD},  Rudelson and Vershynin have proved  that there is a universal constant $c$ \st   the smallest singular value $s_{\min}(z-\bA)$ of $z-\bA$ has order at least $n^{-c}$ as $z$ varies in $\C$ (and stays bounded away from $0$ if $\bA$ is not invertible): this  strong result seems incomplete, as it  does not exhibit any  transition as $|z|$ gets larger than $b$ (by the Single Ring Theorem, we would expect a transition from the order $n^{-c}$ to the order $1$ as $|z|$ gets larger than $b$).  Moreover, one cannot expect  the methods of \cite{RUD}   to allow to prove such a transition, as they are based on the formula $$s_{\min}(z-\bA)\ge \ff{\sqrt{n}}\, \ti\,\min_{1\le i\le n}\op{dist}_2(\trm{ $i$th row of $z-\bA$}, \, \op{span}(\trm{other rows   of $z-\bA$})\,),$$whose RHT cannot have order larger than $n^{-1/2}$.

In this text, we want to fill in the gap of understanding why the borders of the support of $\mu$ have radiuses $b$ and $a$ (it suffices to understand the radius $b$, as $a$ appears then naturally by considering $\bA^{-1}$ instead of $\bA$). For this purpose, by an elementary moment expansion, we shed light on Formula \eqre{9914} by proving that for $k\sim n^{1/6}$, the operator norm $\| \bA^k\|$ has order at most $b^k$. More precisely, in Theorem \re{th8914}, we show that  \be\la{991013h}\E\Tr (\bA^{k}(\bA^{k})^*)\ \lesssim \ b^{2k}.\ee This estimate  allows to state some  exponential bounds for the convergence of the modulii of the extreme eigenvalues of $\bA$ to $a$ and $b$ (see Corollary \re{corSSRT}). As we said above, such a convergence had already been proved  Guionnet and Zeitouni in \cite{GUI2} with some bounds of the type $n^{-\al}$ for an unspecified $\al>0$, but in several applications (as   the study of the outliers related to this matrix model in \cite{FloJean}), polynomial bounds are not enough, while exponential bounds are. We also slightly weaken the hypothesis for this convergence, using the paper \cite{BD13} by Basak and  Dembo.
Then, the estimate \eqre{991013h} allows to give an upper-bound on the rate of convergence of the spectral radius  $|\lam_{\max}(\bA)|$ of $\bA$  to $b$ as $n$ tends to infinity:   Corollary \re{corSSRT1} states  that $$|\lam_{\max}(\bA)|-b \ \lesssim \  n^{-1/6}\log n.$$This result can be compared to the result of Rider in \cite{Rider2003}, who proved that the spectral radius of a Ginibre matrix fluctuates around its limit at rate $(n\log n)^{-1/2}$ (see also some generalizations in \cite{BourgadeCirc2,yinLCL3,DjalilPeche,RiderSinclair}).
At last, in Theorem \re{th8914}, we also prove  \be\la{991013h1}\E[|\Tr (\bA^{k})|^2] \ \lesssim \ b^{2k}\ee with very little efforts, as the proof is mostly  analogous to the one of \eqre{991013h}. The estimate \eqre{991013h1}   is not needed to prove Corollaries  \re{corSSRT} and \re{corSSRT1}, but   will be of use in a forthcoming paper.

The main tools of the proofs are the so-called Weingarten calculus, an integration method for the Haar measure on the unitary group developed by Collins and \'Sniady in \cite{collinsIMRN,collins-sniady06}, together with an exact formula for the Weingarten function (see \eqre{wg2914}) proved by Mastomoto and Novak  in their study of the relation between the Weingarten function and Jucys-Murphy  elements in \cite{MastuNovakIMRN,NovakBanach}. A particularity of this paper is that   the Weingarten calculus is used here to consider products    of Haar-distributed unitary matrices entries with number of factors tending to infinity as the dimension $n$ tends to infinity. \\ \\ 

\section{Main results}
Let $n\ge 1$,  $\bA=\bU\bT\bV$   with $\bT=\diag(s_1, \ld, s_n)$  deterministic \st for all $i$, $s_i\ge 0$,  and $\bU,\bV$  some independent $n\ti n$ Haar-distributed unitary matrices. Set $$  M:=\max_{1\le i\le n}  s_i \qquad\trm{ and }\qquad   b^2:= \ff{n}\sum_{i=1}^n s_i^2 \,.$$

\beg{Th}\la{th8914}Let $\eps>0$. There is a finite constant $C$ depending  only on $\eps$ (in particular, independent of $n$ and of the $s_i$'s) \st for   all positive integer $k$ \st $k^6<(2-\eps)n$, we have
\be\la{59141}\E\Tr (\bA^{k}(\bA^{k})^*) \ \le \ Cnk^2 \lf(b^2+\f{kM^2}{n}\ri)^{k}   \ee
and
\be\la{591412}\E[|\Tr (\bA^{k})|^2] \ \le \ C \lf(b^2+\f{kM^2}{n}\ri)^{k} .  \ee 
\en{Th}

Let now $\lam_{\max}(\bA)$ and  $\lam_{\min}(\bA)$ denote some eigenvalues of $\bA$ with  respectively    largest and smallest    absolute values. We also introduce $\ds a>0$ defined by $$\ds \ff{a^2}:=\ff{n}\sum_{i=1}^n s_i^{-2}$$ (with the convention $\ds \ff{0}=\infty$ and $\ds \ff{\infty}=0$). At last, for $\bM$ a matrix, $\|\bM\|$ denotes the operator norm of $\bM$  with respect to the canonical Hermitian norm.
\beg{cor}\la{corSSRT}With the above notation, there are some constants  $C,\del_0>0$ depending only on $M$ and $b$ \st for any $\del\in [0, \del_0]$,  \be\la{59143}\p(|\lam_{\max}(\bA)|>b+\del)\le Cn^{4/3} \exp\lf(- \f{n^{1/6}\del }{C} \ri)\ee and analogously, if $a>0$ and $m:=\min_i s_i>0$, there are some constants  $C,\del_0>0$ depending only on $a$ and $m$ \st for any $\del\in [0, \del_0]$,  \be\la{59144}\p(|\lam_{\min}(\bA)|<a-\del)\le Cn^{4/3} \exp\lf(- \f{n^{1/6}\del }{C} \ri).\ee
%for any $\eta>0$, there is a constant  $C>0$   \st \be\la{59144}\p(|\lam_{\min}(\bA)|<a-\eta)\le C e^{-\f{n^{1/6}}{C}}.\ee
\en{cor}

\bpr To prove \eqre{59143}, it suffices to notice that if $|\lam_{\max}(\bA)|>b+\del$, then  $|\lam_{\max}(\bA^k)|>(b+\del)^k$, which implies that $\|\bA^k\|>(b+\del)^k$ and that $\Tr (\bA^{k}(\bA^{k})^*)>(b+\del)^{2k}$. Then the Tchebichev inequality and \eqre{59141} allow to conclude.  The proof of \eqre{59144} then follows by application of \eqre{59143} to $\bA^{-1}$ (the matrix $\bA$ is invertible as soon as $a>0$).% (when the $s_i$'s are all positive, as in the other case there is nothing to prove).
%We have \beq  \p(|\lam_{\max}(\bA)|>b+\eta)&\le& \p(|\lam_{\max}(\bA^k)|>(b+\eta)^k)\\
%&\le & \p(\|\bA^k\|>(b+\eta)^k)\\
%&\le & \p(\Tr (\bA^{k}(\bA^{k})^*)>(b+\eta)^{2k})\\
%&\le & \f{\E\Tr (\bA^{k}(\bA^{k})^*)}{(b+\eta)^{2k}}
%\eeq
%and then the theorem allows to conclude.
\epr

The previous results are \emph{non asymptotic} in the sense of \cite{RUDAMSSC,RUDICM}, meaning that they are true for all $n$ (even though they involve some non specified constants). Let us now give  two asymptotic corollaries. The  proof of the first one follows directly from the previous lemma.  Here, $u_n\gg v_n$ means $u_n/v_n\lto 0$.

\beg{cor}\la{corSSRT1}Let now the matrix $\bA=\bA_n$ depend on $n$ and suppose that as $n$ tends to infinity, the numbers $M=M_n$ and $b=b_n$ introduced in Theorem \re{th8914} stay bounded away from $0$ and $+\infty$. Then for any sequence $\del_n>0$,   $$\del_n \, \gg \, n^{-1/6}\log n \; \implies \;  \p(|\lam_{\max}(\bA_n)|>b_n+\del_n) \, \ninf \, 0.$$ The analogous result is also true for  $\lam_{\min}(\bA_n)$.
\en{cor}

 Our last corollary allows a weakenning of the hypotheses of the Single Ring Theorem support convergence proved in \cite{GUI2}  by Guionnet and Zeitouni where we do not even have any single ring anymore. 
 %We set $$\qquad G_\bT(z)=\ff{n}\Tr \ff{z-\bT}\qquad (z\in \C\bck\R).$$
 Let now the matrix $\bA=\bA_n$ depend on $n$, $\bT=\bT_n$ be random, independent from $\bU_n$ and $\bV_n$ and suppose that there is a (possibly random)  
 $M_\infty>1$ independent of $n$ \st with \pro tending to one, the spectrum  of $\bT$ is contained  in $[M_\infty^{-1}, M_\infty]$ and that 
 there is a (possibly random) closed set $K\subset \R$ of zero Lebesgue measure such that for every
$\eps > 0$,  there are some (possibly random) $\ka_\eps > 0$, $M_\eps$   and all $n$ large enough,\be\la{910141}
\{z \in \C \ste \im(z) > n^{-\ka_\eps},\im(\Tr((\bT-z)^{-1}) > n M_\eps\} \subset  \cup_{x\in K} B(x,\eps).\ee
 
\beg{cor}\la{corSSRT91014} $\bullet$ If there   is a finite (possibly random) number   $b\ge 0$ \st as $n\to\infty$,  we have the convergence  in probability $$ \ff{n}\Tr (\bT^2)\lto b^2,$$
then the spectral radius of $\bA$ converges in \pro to $b$.  

$\bullet$ If there is a finite (possibly random) number   $a>0$ \st as $n\to\infty$,  we have the convergences in probability $$  \ff{n}\Tr (\bT^{-2})\lto a^{-2},$$
then the minimal absolute value of the eigenvalues of $\bA$ converges in \pro to $a$. 
 \end{cor}

\bpr
First of all, we shall only prove the first part of the corollary, the proof of the second one being an analogous consequence of Proposition 1.3  of \cite{BD13}. Secondly, up to the replacement of $\bT$ by e.g. $\f{M_\infty+M_\infty^{-1}}{2}I_n$ when its spectrum is not contained  in $[M_\infty^{-1}, M_\infty]$ and to the conditioning with respect to the $\si$-algebra generated by the sequence $\{\bT_n\ste n\ge 1\}$, one can suppose that $\bT$ is deterministic, as well as $M_\infty, b, K, \ka_\eps, M_\eps$, that $\|\bT^{-1}\|,\|\bT\|\le M_\infty$ and that $$\ff{n}\Tr (\bT^2)\lto b.$$ Then as the set of \pro measures supported by $[M_\infty^{-1}, M_\infty] $ is compact, up to an extraction, one can suppose that there is a \pro measure $\Theta$ on $ [M_\infty^{-1}, M_\infty] $ \st the empirical spectral law of $\bT$ converges to $\Theta$ as $n\to\infty$. 
It follows\footnote{In Proposition 1.3  of \cite{BD13} there is the supplementary hypothesis that $\Theta$ is not a Dirac mass, but this restriction might be   there only for the harmonic analysis   characterization of the limit measure true. Indeed,  if $\Theta=\del_b$, then   the convergence of the empirical spectral law of $\bA$ to the uniform measure on the circle with radius $b$ is obvious from the convergence of the empiricial spectral distribution of a Haar-distributed unitary matrix to the uniform law on the unit circle and from classical perturbation inequalities.}, by Proposition 1.3  of \cite{BD13}, that within a subsequence, the empirical spectral law of $\bA$ converges in \pro to a \pro measure on $\C$ whose support is a single ring, with maximal radius
$\int x^2\ud \Theta(x)=b^2$, hence that for any $\eps>0$, $$\p(|\lam_{\max}(\bA)|<b-\eps)\lto 0.$$ By \eqre{59143}, the convergence in \pro of the spectral radius to $b$ is proved within a subsequence. In fact, we have even proved more: we have proved that from any subsequence, we can extract a subsequence within which the spectral radius of $\bA$ converges in \pro to $b$. This is enough to conclude.  \epr

\section{Proof of Theorem \re{th8914}} We  first prove \eqre{59141} (we will see below that the proof of \eqre{591412} will go along the same lines, minus some   border  difficulties).

We have\beq
\E\Tr (\bA^{k}(\bA^{k})^*)&=&\E\Tr \underbrace{\bU\bT\bV\cd \bU\bT\bV}_{\trm{$k$ times}}\underbrace{\bV^*\bT\bU^*\cd \bV^*\bT\bU^*}_{\trm{$k$ times}}\\
&=&\E\Tr \underbrace{\bT\bV\bU\cd \bT\bV\bU}_{\trm{$k-1$ times}}\bT^2\underbrace{(\bV\bU)^*\bT\cd (\bV\bU)^*\bT}_{\trm{$k-1$ times}}\\
&=&\E\Tr \underbrace{\bT\bU\cd \bT\bU}_{\trm{$k-1$ times}}\bT^2\underbrace{\bU^*\bT\cd \bU^*\bT}_{\trm{$k-1$ times}}\, ,\eeq where we used  the fact that $\bV\bU\eqlaw \bU$.

Let us denote $\bU=[u_{ij}]_{i,j=1}^n$
 and $\bU^*=[u^*_{ij}]_{i,j=1}^n$. Then continuing the previous computation, we get \beq&&
\E\Tr (\bA^{k}(\bA^{k})^*)=\\ &&\sum_{\substack{\bi=(i_1, \ld, i_{k })\\ \bj=(j_1, \ld, j_{k })\\ j_1=i_1, j_{k }=i_{k }}} \E s_{i_1}u_{i_1i_2}s_{i_2}u_{i_2i_3}\cd s_{i_{k-1}}u_{i_{k-1}i_{k}}s_{i_{k}}s_{j_{k}}u^*_{j_{k}j_{k-1}}s_{j_{k-1}}u^*_{j_{k-1}j_{k-2}}s_{j_{k-2}}\cd u^*_{j_{2}j_{1}}s_{j_{1}}\,.\eeq By left and right invariance of the Haar measure (see Proposition \re{wg}), for the expectation in the RHT to be non zero, we need to have the equality of multisets $$ \{j_1, \ld, j_{k}\}_m=\{i_1, \ld, i_{k}\}_m$$ (the subscript $m$ is used to denote multisets here).
So \beq
&& \!\!\!\!\!\!\!\!\!\!\E\Tr (\bA^{k}(\bA^{k})^*)
=\\ \\ &&\sum_{\bi=(i_1, \ld, i_{k })}\!\!\!\!\!\sum_{\substack{\bj=(j_1, \ld, j_{k })\\ j_1=i_1, j_{k }=i_{k }\\ \{j_2, \ld, j_{k-1}\}_m=\{i_2, \ld, i_{k-1}\}_m}}\!\!\!\!\!\!\!\!\!\!\!\!\!\!\!\!\!\!\!\!s_{i_1 }\cd s_{i_{k }}s_{j_{k }}\cd s_{j_{1}} \E u_{i_1i_2}u_{i_2i_3}\cd  u_{i_{k-1}i_{k }} u^*_{j_{k }j_{k-1}}s_{j_{k-1}}u^*_{j_{k-1}j_{k-2}} \cd u^*_{j_{2}j_{1}}\\ &&\\ &&=\sum_{\bi=(i_1, \ld, i_{k })}\!\!\!\!\!\!\!\!\!\sum_{\substack{\bj=(j_1, \ld, j_{k })\\ j_1=i_1, j_{k }=i_{k }\\ \{j_2, \ld, j_{k-1}\}_m=\{i_2, \ld, i_{k-1}\}_m}} \!\!\!\!\!\!\!\!\!s_{i_1}^2\cd s_{i_{k}}^2\E  u_{i_1i_2} u_{i_2i_3}\cd  u_{i_{k-1}i_{k }} u^*_{j_{k }j_{k-1}} u^*_{j_{k-1}j_{k-2}} \cd u^*_{j_{2}j_{1}} \,.
\eeq	
Let us define the subgroup $\Sy_{k}^0$ of the $k$th symmetric group $\Sy_{k}$ by \beq \Sy_{k}^0&:=&\{\vfi\in \Sy_{k}\ste \vfi(1)=1, \, \vfi(k )=k \},\eeq
and for each $\bi=(i_1, \ld, i_{k})$, we define the stabilisator group \beq\Sy_{k}^0(\bi)&:=&\{\al\in \Sy_{k}^0\ste \;\forall \ell, \;i_{\al(\ell)}=i_\ell\}\,.\eeq
Let  at last  $\Sy_{k}^{0}/\Sy_{k}^{0}(\bi)$ denote the quotient of the set $\Sy_{k}^{0}$ by $\Sy_{k}^{0}(\bi)$ for the left 
action $$(\al, \vfi)\in \Sy_{k}^{0}(\bi)\ti\Sy_{k}^{0}\longmapsto \al\vfi.$$
Remark  that the notation $i_{\Phi(1)}, \ld, i_{\Phi(k)}$ makes sense for $\Phi\in   \Sy_{k}^{0}/\Sy_{k}^{0}(\bi)$ even though $\Phi$ is not a permutation but a set of permutations. 
Then we have 
 \be\la{89141h}
\E\Tr (\bA^{k}(\bA^{k})^*) = \sum_{\bi=(i_1, \ld, i_{k })\in \{1,  \ld, n\}^{k }}(s_{i_1}^2\cd s_{i_{k}}^2\,\ti\,F_\bi)
 \ee
with   \be\la{591417h19}F_\bi:= \sum_{\Phi\in \Sy_{k}^{0}/\Sy_{k}^{0}(\bi)} \E  u_{i_1i_2} u_{i_2i_3}\cd  u_{i_{k-1}i_{k }} u^*_{i_{\Phi(k)}i_{\Phi(k-1)}} \cd u^*_{i_{\Phi(2)}i_{\Phi(1)}} \,.
  \ee
%We obviously have $$\E  |u_{i_1i_2}|^2 |u_{i_2i_3}|^2\cd  |u_{i_{k-1}i_k}|^2=\wg(id)$$ and i
%Note that for any $\Phi\in \Sy_{k}^0$, we have $ \Phi c^{-1} \Phi^{-1}c(k+1)=k+1,$ so that $$ \Phi c^{-1} \Phi^{-1}c$$  defines a unique element $f(\Phi)\in \Sy_{k}$.
Let us apply Proposition \re{wg} to compute the expectation in the term associated to an element  $\Phi\in \Sy_{k}^{0}$. Let $\wg$ denote the Weingarten function, introduced in Proposition \re{wg}. Let $c\in \Sy_{k}$ be the cycle $(12\cd k)$. %Note that for any $\tau\in \Sy_k$, $c\tau c^{-1}$ defines a permutation of the set $\{2, \ld, k+1\}$, that $c^{-1}\Phi^{-1}c$ defines a permutation of the set $\{1, \ld, k\}$ and that $$ \Phi c^{-1} \Phi^{-1}c$$  defines a unique element $f(\Phi)\in \Sy_{k}$.
For any $\bi=(i_1, \ld, i_{k })$ and any $\Phi\in \Sy_{k}^{0}$,%/\Sy_{k}^0(\bi)$, 
 \beq &&\E  u_{i_1i_2} u_{i_2i_3}\cd  u_{i_{k-1}i_{k }} u^*_{i_{\Phi(k)}i_{\Phi(k-1)}} \cd u^*_{i_{\Phi(2)}i_{\Phi(1)}}\ =\ 
\\
&&\qquad\qquad\qquad\qquad  \E  u_{i_1i_2} u_{i_2i_3}\cd  u_{i_{k-1}i_{k }}  \ovl{u_{i_{\Phi(1)}i_{\Phi(2)}}}\cd \ovl{u_{i_{\Phi(k-1)}i_{\Phi(k)}}}\ =\ 
\\
&&\qquad\qquad\qquad\qquad   \sum_{\sigma, \tau \in \Sy_{k-1}} \delta_{i_1,i_{\Phi\sigma(1)}} \ldots \delta_{i_{k-1},i_{\Phi\sigma(k-1)}} \delta_{i_2,i_{\Phi c\tau c^{-1}(2)}} \ldots \delta_{i_{k },i_{\Phi c\tau c^{-1}(k )}} \wg( \si^{-1}\tau)\,.\eeq
But for any $\Phi\in \Sy_{k}^{0}$ and any $\si,\tau\in \Sy_{k-1}$, with $\ovl{\si},\ovl{\tau}\in \Sy_{k}$ defined by 
$$\ovl{\si}(x)=\beg{cases} \si(x)&\trm{ if $1\le x\le k-1$, }\\ k &\trm{ if $x=k $, }\en{cases}\qquad\qquad \ovl{\tau}(x)=\beg{cases} \tau(x)&\trm{ if $1\le x\le k-1$, }\\ k &\trm{ if $x=k $, }\en{cases}$$
  we have (with the conventions that    $(x\, y)$ denotes the transposition of $x$ and $y$ when $x\ne y$ and the identity otherwise and that for $\pi\in \Sy_k$ \st $\pi(k)=k$, $\pi_{|_{\{1, \ld, k-1\}}}$ denotes the restriction of $\pi$ to $\{1, \ld, k-1\}$)
\beq  i_1=i_{\Phi\sigma(1)},\ldots,  i_{k-1}=i_{\Phi\sigma(k-1)}&\iff&  i_1=i_{\Phi\ovl{\sigma}(1)},\ldots,  i_{k}=i_{\Phi\ovl{\sigma}(k)}\\
&\iff&\trm{for }\ell_1=\si^{-1}(1),\; \forall \vfi\in \Phi, \; \vfi\ovl{\si}\,(1\,\ell_1)\in 
  \Sy_{k}^0(\bi)\\ &\iff &\exists \ell_1\in \{1, \ld, k-1\},\; \forall \vfi\in \Phi, \; \vfi\ovl{\si}\, (1\,\ell_1)\in 
  \Sy_{k}^{0}(\bi)
  \\ &\iff & \exists \ell_1\in \{1, \ld, k-1\},\;\exists \vfi_1\in \Phi, \;\si=(\vfi_1^{-1}(1\,\ell_1))_{|_{\{1, \ld, k-1\}}} 
\eeq
and 
\beq i_2=i_{\Phi c\tau c^{-1}(2)}, \ldots ,i_{k }=i_{\Phi c\tau c^{-1}(k )}&\iff&i_1=i_{\Phi c\ovl{\tau} c^{-1}(1)},  \ldots ,i_{k }=i_{\Phi c\ovl{\tau} c^{-1}(k )}\\ &\iff&\trm{ for }\ell_2=c\tau^{-1}(k-1),\;  \forall \vfi\in \Phi, \\&&\qquad\qquad\vfi c\ovl{\tau} c^{-1}\,(\ell_2\, k)\in \Sy_{k}^0(\bi)\\ &\iff &\exists \ell_2\in \{2, \ld, k \},\; \forall \vfi\in \Phi, \\&&\qquad\qquad \vfi c\ovl{\tau} c^{-1}(\ell_2\,  k \,)\in \Sy_{k}^{0}(\bi)
\\ &\iff &\exists \ell_2\in \{2, \ld, k \},\;  \exists \vfi_2\in \Phi, \\&&\qquad\qquad c\ovl{\tau} c^{-1}=\vfi_2^{-1} (\ell_2\,  k \,)\\
&\iff &\exists \ell_2\in \{1, \ld, k-1\},\;  \exists \vfi_2\in \Phi, \\&&\qquad\qquad \ovl{\tau} =c^{-1}\vfi_2^{-1}c\,  (\ell_2\, k-1)
\\ &\iff & \exists \ell_2\in \{1, \ld, k-1\},\;\exists \vfi_2\in \Phi, \\ &&\qquad\qquad  \tau =(c^{-1}\vfi_2^{-1}c\,  (\ell_2\, k-1) )_{|_{\{1, \ld, k-1\}}} .
\eeq
Hence %\beq \E  u_{i_1i_2} u_{i_2i_3}\cd  u_{i_{k}i_{k+1}} u^*_{i_\Phi(k)i_{\Phi(k)}} \cd u^*_{i_{\Phi(2)}i_{\Phi(1)}}&=&\sum_{(\si,\tau)\in \Phi^{-1}\Sy_{k}^0(\bi)\ti c^{-1}\Phi^{-1}\Sy_{k}^0(\bi)c}\wg( \si^{-1}\tau)\\ &=&\sum_{h,h'\in \Sy_{k}^0(\bi)}\wg(h^{-1}\Phi c^{-1}\Phi^{-1}h'c)\eeq
 \beq &&\E  u_{i_1i_2} u_{i_2i_3}\cd  u_{i_{k}i_{k+1}} u^*_{i_\Phi(k)i_{\Phi(k)}} \cd u^*_{i_{\Phi(2)}i_{\Phi(1)}}=
\\ &&\qquad\qquad\qquad\qquad\sum_{\substack{(\vfi_1, \vfi_2)\in \Phi\ti\Phi\\ 1\le \ell_1, \ell_2\le k-1}}\wg( ((1\,\ell_1)\vfi_1c^{-1}\vfi_2^{-1}c\, (\ell_2\, k-1))_{|_{\{1, \ld, k-1\}}} )=\\
&&\qquad\qquad\qquad\qquad\sum_{\substack{(\vfi_1, \vfi_2)\in \Phi\ti\Phi\\ 1\le \ell_1, \ell_2\le k-1}}\wg((c^{-1}\vfi_2^{-1}c\, (\ell_2\, k-1)(1\,\ell_1)\vfi_1)_{|_{\{1, \ld, k-1\}}})\,, \eeq
where we used the fact that $\wg$ is a central function. 

Thus by the definition of $F_\bi$ at \eqre{591417h19},  we have  \beqy \nonumber F_\bi 
&=& \sum_{\Phi\in \Sy_{k}^{0}/\Sy_{k}^{0}(\bi)}   \sum_{\substack{(\vfi_1, \vfi_2)\in \Phi\ti\Phi\\ 1\le \ell_1, \ell_2\le k-1}}\wg((c^{-1}\vfi_2^{-1}c\, (\ell_2\, k-1)(1\,\ell_1)\vfi_1)_{|_{\{1, \ld, k-1\}}})\\
\nonumber&=& \sum_{1\le \ell_1,\ell_2\le k-1}\sum_{\substack{(\vfi_1, \vfi_2)\in \Sy_{k}^{0}\ti\Sy_{k}^{0}\\ \vfi_1=\vfi_2\trm{ in } \Sy_{k}^{0}/\Sy_{k}^{0}(\bi)}}    \wg( (c^{-1}\vfi_2^{-1}c\, (\ell_2\, k-1)(1\,\ell_1)\vfi_1)_{|_{\{1, \ld, k-1\}}} )\\
\la{591417h18}&=&  \sum_{1\le \ell_1,\ell_2\le k-1}\sum_{(\vfi, \al) \in \Sy_{k}^{0}\ti\Sy_{k}^{0}(\bi) }    \wg( (c^{-1}\vfi^{-1}\al^{-1}c\, (\ell_2\, k-1)(1\,\ell_1)  \vfi)_{|_{\{1, \ld, k-1\}}} )\,.
  \eeqy

  To state the following lemma, we need to introduce some notation: for  $\si\in \Sy_k$, let $|\sigma|$ denote the minimal number of factors necessary to write $\sigma$ as a product of transpositions.
  \beg{lem}\la{majwg} Suppose that  $k^2<2n$. Then for any $\pi\in \Sy_k\bck\{id\}$, we have \be\la{591413h} |\wg(\pi)|\le \f{2}{n^kk^2}\lf(\f{k^2}{2n}\ri)^{|\pi|}\ff{\ds 1-\f{k^2}{2n}}\ee and  \be\la{591413h1} |\wg(id)|\le\ff{n^k}+\f{k^2}{2n^{k+2}} \ff{\ds 1-\f{k^4}{4n^2}}\,.\ee
  \end{lem}

\bpr
We know, by \eqre{wg2914}, that the (implicitly depending on $n$) function $\wg$ can be written $$\wg(\pi)=\ff{n^k}\sum_{r\ge 0}(-1)^r\f{c_r(\pi)}{n^r},$$ with $$c_r(\pi):=\#\{(s_1, \ld, s_r, t_1, \ld, t_r)\ste \forall i, 1\le s_i<t_i\le k,\;t_1\le\cd\le t_r,\;\pi=(s_1\,t_1)\cd (s_r\,t_r)\}.$$
But  for $r\ge 1$,  $$c_r(\pi)\le \one_{r\ge |\pi|}\binom{k}{2}^{r-1}\le  \one_{r\ge |\pi|}\f{k^{2r-2}}{2^{r-1}},$$ 
so  that for $\pi\ne id$,  $$ |\wg(\pi)| \ \le \ \f{2}{n^kk^2}\sum_{r\ge |\pi|}\lf(\f{k^2}{2n}\ri)^r \ \le \  \f{2}{n^kk^2}\lf(\f{k^2}{2n}\ri)^{|\pi|}\ff{\ds 1-\f{k^2}{2n}},
$$
whereas 
$$|\wg(id)| \ \le  \ \ff{n^k}+\ff{n^k}\sum_{r\ge 1} \f{|c_{2r}(id)|}{n^{2r}} \ \le \ \ff{n^k}+\f{2}{n^kk^2}\sum_{r\ge 1}\lf(\f{k^2}{2n}\ri)^{2r} \ = \  \ff{n^k}+\f{k^2}{2n^{k+2}} \ff{\ds 1-\f{k^4}{4n^2}}.$$
\epr

\beg{rem}Some other upper-bounds have been given for the Weingarten function: Theorem 4.1 of \cite{CollinsSpain} and 
  Lemma 16 of  \cite{Montanaro}. The first one states that for any $k,j\ge 2$  \st $k^j\le n$, there is $K_j$ depending only on $j$ \st for   all $\pi\in   \Sy_k$, \be\la{Collinsandco29914}|\wg(\pi)|\le K_jn^{-k-|\pi|(1-2/j)},\ee whereas the second one states that if $k^{3/2}\le n$, then  for     all $\pi\in   \Sy_k$, \be\la{Collinsandco299142}|\wg(\pi)|\le \f{3C_{k-1}}{2}n^{-k-|\pi|},\ee
where $C_{k-1}$ is the Catalan number of index $k-1$. However, in our case, these bounds are less relevant than the one we give here. Indeed,   \eqre{Collinsandco299142} allows to weaken the hypothesis $k^6\le (2-\eps)n$ to $k^4\le (2-\eps)n$, but, because of the Catalan number, contains implicitely a factor $4^k$, which would change our main result from $$\E\Tr (\bA^{k}(\bA^{k})^*)   \ \lesssim \  b^{2k}$$ to $$\E\Tr (\bA^{k}(\bA^{k})^*)   \ \lesssim \  (2b)^{2k}$$ (which is far less interesting in our point of view, as explained in the introduction).
On the other hand, \eqre{Collinsandco29914} allows to turn the hypothesis $k^6\le (2-\eps)n$ to $$k^{\max\{j,\, 4j/(j-2)\}}\le (2-\eps)n,$$  but there is no integer $j$ making it a good deal. 
\en{rem}

It follows from Lemma \re{majwg}  (that we apply for   $k-1$ instead of $k$) and from \eqre{591417h18}  that  \beqy\nonumber |F_\bi|&\le &  \sum_{1\le \ell_1,\ell_2\le k-1}\sum_{\al\in \Sy_{k}^{0}(\bi)} \Bigg\{\\&&\nonumber\qquad  \#\{\vfi\in \Sy_{k}^{0}\ste  c^{-1}\vfi^{-1}\al^{-1}c\, (\ell_2\, k-1)(1\,\ell_1)  \vfi =id\}\bigg(\ff{n^{k-1}}+\f{k^2}{2n^{k+1}} \ff{\ds 1-\f{k^4}{4n^2}}\bigg)+\\ \la{591417h191}&&   \ff{\ds 1-(k^2/(2n))}\sum_{q=1}^{k-2}\bigg( \#\{\vfi\in \Sy_{k}^{0}\ste |c^{-1}\vfi^{-1}\al^{-1}c\, (\ell_2\, k-1)(1\,\ell_1)  \vfi|=q\} \ti\nonumber \\ &&  \qquad\qquad\qquad\qquad\qquad \qquad\qquad\qquad \qquad\qquad\qquad\qquad \f{2}{n^{k-1}k^2}\lf(\f{k^2}{2n}\ri)^{q}\bigg) \Bigg\} \eeqy (note that %$c^{-1}\vfi^{-1}\al^{-1}c\, (\ell_2\, k-1)(1\,\ell_1)  \vfi$ always fixes $k $, this is why we can 
we suppress the restriction to the set $\{1,\ld, k-1\}$ because adding a fixed point does not change the minimal number of transposition needed to write a permutation).

Thus to upper-bound $|F_\bi|$, we need to upper-bound, for $1\le \ell_1,\ell_2\le k-1$, $\al\in \Sy_{k}^{0}(\bi)$ and $q\in \{0,\ld, k-2\}$ fixed, the cardinality of the set of $\vfi$'s in $\Sy_{k}^{0}$ \st $$ |c^{-1}\vfi^{-1}\al^{-1}c\, (\ell_2\, k-1)(1\,\ell_1)  \vfi|=q.$$

\beg{lem}\la{majwg2}
Let $q\in \{0, \ld, k-2\}$, $1\le \ell_1,\ell_2\le k-1$ and $\al\in \Sy_{k}^{0}$ be fixed.  Then we have 
\be\la{591416h}\#\{\vfi\in \Sy_{k}^{0}\ste |c^{-1}\vfi^{-1}\al^{-1}c\, (\ell_2\, k-1)(1\,\ell_1)  \vfi|=q\} \ \le \ \f{k^{4q}}{(2q)!}.\ee\en{lem}

\bpr Let $\vfi\in\Sy_k^0$ and define $$\pi_\vfi:=c^{-1}\vfi^{-1}\al^{-1}c\, (\ell_2\, k-1)(1\,\ell_1)   \vfi.$$ Note that for any $x\in \{1, \ld, k \}$, $$
\pi_\vfi(x)=x \ \implies \ \vfi(x)=(1\,\ell_1) (\ell_2\, k-1)c^{-1}\al\vfi c(x),$$ so that $\vfi(x)$ is determined by $\vfi(c(x))$ (and by $\ell_1$, $\ell_2$ and $\al$, which are considered as fixed here). 

$\bullet $ It first follows that if $|\pi_\vfi|=0$, \ie if all $x$'s are fixed points of $\pi_\vfi$, then as by definition of $\Sy_{k}^{0}$, we always have $\vfi(k )=k $, $\vfi$ is entirely defined by $\al$, $\ell_1$ and $\ell_2$. It proves \eqre{591416h} for $q=0$. 

$\bullet$ (i) If $|\pi_\vfi|\ne 0$, \ie if $\pi_\vfi$ has not only fixed points, then it follows that the values of $\vfi$ on the complementary of the support of $\pi_\vfi$ are entirely determined by its values on the support of $\pi_\vfi$ (and by $\ell_1$, $\ell_2$ and $\al$). 

(ii) Let us now define, for each $\si\in \Sy_{k}$ \st $|\si|=q$,  a subset $\A(\vfi)$ of $\{1,\ld, k\}$   \st $$  \#\A(\si) =  2q \qquad\trm{ and }\qquad\supp(\si)\subset \A(\si).$$
Such a set can be defined as follows: any permutation $\si$ at distance $q$ from the identity admits one (and only one, in fact, but we do not need this here) factorization 
	$$\si = (s_1 t_1) \cdots (s_q t_q)$$
such that $s_i < t_i$ and $t_1 < \dots < t_q$. One can choose $\A(\si)$ to be $\{s_1, t_1, \ldots ,s_q, t_q\}$, possibly arbitrarily completed to a set with cardinality $2q$ if needed.

(iii) By what precedes,  the map \beq\{\vfi\in \Sy_{k}^{0}\ste |\pi_\vfi|=q\}&\lto&\underset{\substack{A\subset\{1,\ld, k\}\\ \#A=2q}}{\cup}\{1, \ld, k\}^A
 \\
\vfi&\longmapsto& \vfi_{|_{\A\lf( \pi_\vfi \ri)}}
\eeq is one-to-one, and 
$$\#\{\vfi\in \Sy_{k}^{0}\ste |\pi_\vfi|=q\}\le  \binom{k}{2q}k^{2q}\le  \f{k^{4q}}{(2q)!}.$$
\epr

It follows from this lemma and from \eqre{591417h191} that for a constant $C$ depending only on the $\eps$ of the statement of the theorem (this constant might change from line to line),  $$ |F_\bi| \ \le \ k^2\#\Sy_{k}^0(\bi) \Bigg\{  \ff{n^{k-1}}+\f{Ck^2}{ n^{k+1}}  +\f{C}{n^{k-1}k^2} \sum_{q=1}^{k-1} \ff{(2q)!} \lf(\f{k^6}{2n}\ri)^{q} \Bigg\} \ \le \  \f{Ck^2}{n^{k-1}} \#\Sy_{k}^{0}(\bi).$$

Thus by \eqre{89141h}, 
$$\E\Tr (\bA^{k}(\bA^{k})^*) \  \le \  \f{Ck^2}{n^{k-1}} \sum_{\bi=(i_1, \ld, i_{k })\in \{1,  \ld, n\}^{k }}(s_{i_1}^2\cd s_{i_{k}}^2\, \ti\,  \#\Sy_{k}^{0}(\bi)).$$
Let us rewrite the last sum as follows: \bgt\ite  we first choose the number $p\in  \{1, \ld, k \}$ \st $\#\{i_1, \ld, i_{k }\}=p$,
\ite  then we choose  the set $\{i^0_1<\cd<i^0_p\}\subset\{1,  \ld, n\} $ \st we have the equality of sets $\{i_1, \ld, i_{k }\}=\{i^0_1,\ld,i^0_p\}$ (we have $\binom{n}{p}$ possibilities),
\ite then we choose  the collection $(\lam_1, \ld, \lam_p)$ of positive integers summing up to $k $ \st in $(i_1, \ld, i_{k })$, $i^0_1$ appears $\lam_1$ times, \ld,  $i^0_p$ appears $\lam_p$ times (we have $\binom{k-1}{p-1}$ possibilities),
\ite at last we choose a collection $\mathbf{S}=(S_1, \ld, S_p)$ of pairwise disjoint subsets of $\{1, \ld, k \}$ whose union is $\{1, \ld, k \}$ and  with respective cardinalities $\lam_1$, \ld, $\lam_p$ (we have $\f{k!}{\lam_1!\cd \lam_p!}$ possibilities).
\ent
The corresponding collection $\bi=(i_1,   \ld, i_k)$ is then totally defined by the fact that for all $\ell$, $i_\ell=i_{r}^0$, where $r\in \{1, \ld, p\}$ is \st $\ell\in S_r$. 
Note that in this case, $$\#\Sy_{k}^{0}(\bi)\le \lam_1!\cd \lam_p!\,.$$

This gives 
\beq \E\Tr (\bA^{k}(\bA^{k})^*) &\le &\f{Ck^2}{n^k}\sum_{p=1}^{k }\sum_{1\le i^0_1<\cd<i^0_p\le n}\sum_{(\lam_1, \ld, \lam_p)}\sum_{\mathbf{S}}\prod_{\ell=1}^ps_{i_\ell}^{2\lam_\ell}\prod_{\ell=1}^p\lam_\ell !\,.
\eeq
Then, note that $$\prod_{\ell=1}^ps_{i_\ell}^{2\lam_\ell} \ = \ \prod_{\ell=1}^ps_{i_\ell}^{2}s_{i_\ell}^{2\lam_\ell -2} \ \le \ M^{2(k -p)}\prod_{\ell=1}^ps_{i_\ell}^{2}.$$
Hence changing the order of summation, we get 
\beq \E\Tr (\bA^{k}(\bA^{k})^*) &\le &\f{Ck^2}{n^{k-1}}\sum_{p=1}^{k }M^{2(k -p)}\sum_{1\le i^0_1<\cd<i^0_p\le n}\prod_{\ell=1}^ps_{i_\ell}^{2}\sum_{(\lam_1, \ld, \lam_p)}\sum_{\mathbf{S}} \prod_{\ell=1}^p\lam_\ell !\\
&\le & \f{Ck^2}{n^{k-1}}\sum_{p=1}^{k }M^{2(k -p)}\ff{p!}\sum_{1\le i^0_1,\ld,i^0_p\le n}\prod_{\ell=1}^ps_{i_\ell}^{2} \sum_{(\lam_1, \ld, \lam_p)} k !\\
&\le & \f{Ck^2}{n^{k-1}}\sum_{p=1}^{k }M^{2(k -p)}\lf(nb^2\ri)^p\f{k!(k-1)!}{p!(p-1)!(k -p)!}\\
&\le & \f{Ck^2}{n^k}\sum_{p=1}^{k }(kM^{2})^{k -p}\lf(nb^2\ri)^p\f{k!}{p!(k -p)!}\\
&\le & \f{Ck^2}{n^{k-1}}\lf(nb^2+kM^2\ri)^{k }\\
&\le & Cnk^2 \lf(b^2+\f{kM^2}{n}\ri)^{k}\,.
\eeq

 Let us now give the main lines of the proof of \eqre{591412}. This proof is very analogous to the one of \eqre{59141}, minus some  border  difficulties (we will use the symmetric group $\Sy_k$ instead of $\Sy_k^0$). 
 
 Proceeding as above, we arrive easily at 
 \beq \E[|\Tr (\bA^{k})|^2] &=& \sum_{\bi=(i_1, \ld, i_k)}s_{i_1}^2\cd s_{i_k}^2 G_\bi\eeq with $$G_\bi \ := \ \sum_{\Phi\in \Sy_k/\Sy_k(\bi)}  \E[u_{i_1i_2}u_{i_2i_3}\cd u_{i_{k-1}i_i}u_{i_ki_1}\ovl{u_{i_{\Phi(1)}i_{\Phi(2)}}}\ovl{u_{i_{\Phi(2)}i_{\Phi(3)}}}\cd  \ovl{u_{i_{\Phi(k-1)}i_{\Phi(k)}}}\ovl{u_{i_{\Phi(k)}i_{\Phi(1)}}}],$$ where $$\Sy_k(\bi):=\{\vfi\in \Sy_k\ste \;\forall \ell, \; i_{\vfi(\ell)}=i_\ell\}$$ and $  \Sy_k/\Sy_k(\bi)$ denotes the quotient set for the left action of $\Sy_k(\bi)$ on $\Sy_k$.
 
 As above again, we get, for $c$ the cycle $(1\, 2\cd k)$,  \beq G_\bi & = &  \sum_{\Phi\in \Sy_k/\Sy_k(\bi)}\sum_{(\vfi_1, \vfi_2)\in \Phi\ti\Phi}\wg(c^{-1}\vfi_2^{-1}c\vfi_1)\\
 &=&\sum_{(\vfi, \al)\in \Sy_k\ti\Sy_k(\bi)}\wg(c^{-1}\vfi^{-1}\al^{-1}c\vfi)\,.\eeq

Then applying Lemma \re{majwg},  we get 
\beq  |G_\bi|&\le &   \sum_{\al\in \Sy_{k}(\bi)} \Bigg\{  \#\{\vfi\in \Sy_{k}\ste  c^{-1}\vfi^{-1}\al^{-1}c   \vfi =id\}\bigg(\ff{n^{k}}+\f{k^2}{2n^{k+2}} \ff{\ds 1-\f{k^4}{4n^2}}\bigg)+\\  &&   \ff{\ds 1-(k^2/(2n))}\sum_{q=1}^{k-1} \#\{\vfi\in \Sy_{k} \ste |c^{-1}\vfi^{-1}\al^{-1}c  \vfi|=q\} \f{2}{n^{k}k^2}\lf(\f{k^2}{2n}\ri)^{q} \Bigg\} \,.\eeq

Then an analogue of Lemma \re{majwg2} allows to claim that $$|G_\bi|\ \le \ \f{C}{n^k}\#\Sy_k(\bi),$$ and the end of the proof  is quite analogous to what we saw above.
\hfill $\square$\\

\section{Appendix: Weingarten calculus}

We   recall this key-result about integration with respect to the Haar measure on the unitary group (see \cite[Cor. 2.4]{collins-sniady06} and \cite[p. 61]{Novak}). Let $\Sy_k$ denote the $k$th symmetric group. \begin{propo} \label{wg}
\indent Let $k$ be a   positive integer and $\bU=[u_{ij}]_{i,j=1}^n$ a Haar-distributed matrix on the unitary group. Let $\bi=(i_1,\ldots,i_k,i'_1,\ldots,i'_k)$, $\bj=(j_1,\ldots,j_k,j'_1,\ldots,j'_k)$ be two $2k$-uplets of $\left\{1,\ldots,n \right\}$.  Then 
\be \label{wg1}
\E \left[ u_{i_1,j_1}  \cdots u_{i_k,j_k} \ovl{ u_{ i'_1 ,j'_1 } }  \cdots \ovl{u_{i'_k,j'_k}} \right]\ = \ \sum_{\sigma, \tau \in \Sy_k} \delta_{i_1,i'_{\sigma(1)}} \ldots \delta_{i_k,i'_{\sigma(k)}} \delta_{j_1,j'_{\tau(1)}} \ldots \delta_{j_k,j'_{\tau(k)}} \wg( \si^{-1}\tau),
\ee
where $\wg$ is a function called  the \emph{Weingarten function}, depending implicitly on $n$ and $k$ and given by the fomula \be \label{wg2914}\wg(\pi)=\ff{n^k}\sum_{r\ge 0}(-1)^r\f{c_r(\pi)}{n^r},\ee with $$c_r(\pi):=\#\{(s_1, \ld, s_r, t_1, \ld, t_r)\ste \forall i, 1\le s_i<t_i\le k,\;t_1\le\cd\le t_r,\;\pi=(s_1\,t_1)\cd (s_r\,t_r)\}.$$
\end{propo}

 \noindent{\bf Acknowledgments:}   We  would like to  thank J. Novak  for   discussions on Weingarten calculus and Camille Male for pointing  out  references \cite{CollinsSpain} and \cite{Montanaro} to us.

\en{document}
\begin{thebibliography}{10}
  
  \bibitem{BD13} A. Basak, A. Dembo \emph{Limiting spectral distribution of sums of unitary and orthogonal matrices}. Electron. Commun. Probab. 18 (2013), no. 69, 19 pp.
  
  \bibitem{FloJean} F. Benaych-Georges, J. Rochet \emph{Outliers in the Single Ring Theorem}, arXiv:1308.3064. 
  
    \bibitem{BourgadeCirc2} P. Bourgade, H.-T. Yau, J. Yin  \emph{The local circular law II: the edge case},    Probab. Theory Related Fields, Vol. 159, no. 3-4, 619--660 (2014).
  
 \bibitem{DjalilPeche} D. Chafa\"\i, S. P\'ech\'e \emph{A note on the second order universality at the edge of Coulomb gases on the plane}. J. Stat. Phys. 156 (2014), no. 2, 368--383.
  
  \bibitem{collinsIMRN} B. Collins \emph{Moments and cumulants of polynomial random variables on unitary groups, the Itzykson-Zuber integral, and free probability}. Int. Math. Res. Not. 2003, no. 17, 953--982.
  
    \bibitem{CollinsSpain} B. Collins, C. E. Gonz\'alez-Guill\'en, D. P\'erez-Garc\'{\i}a \emph{Matrix Product States, Random Matrix Theory and the Principle of Maximum Entropy}, Comm. Math. Phys. 320 (2013), no. 3, 663--677. 
  
\bibitem{collins-sniady06}  B. Collins, P. \'Sniady \emph{Integration with respect to the Haar measure on unitary, orthogonal and symplectic group}. Comm. Math. Phys., 264, 773--795, 2006.

 \bibitem{GUI} A. Guionnet, M. Krishnapur, O. Zeitouni \emph{The Single Ring Theorem}. Ann. of Math. (2) 174 (2011), no. 2, 1189--1217.

\bibitem{GUI2} A. Guionnet, O. Zeitouni \emph{Support convergence in the Single Ring Theorem}. Probab. Theory Related Fields 154 (2012), no. 3-4, 661--675.

\bibitem{haag2}   U. Haagerup, F.  Larsen \emph{Brown's spectral distribution measure for $R$-diagonal elements in finite von Neumann algebras},  {Journ. Functional Analysis}, {176} (2000), 331--367.

%\bibitem{hiai} F. Hiai, D. Petz \emph{The semicircle law, free random variables, and entropy}. Amer. Math. Soc., Mathematical Surveys and Monographs Volume 77, 2000.

\bibitem{MastuNovakIMRN} S. Matsumoto, J. Novak \emph{Jucys-Murphy elements and unitary matrix integrals}. Int. Math. Res. Not. IMRN 2013, no. 2, 362--397.

\bibitem{Montanaro} A. Montanaro \emph{Weak multiplicativity for random quantum channels},  Comm. Math. Phys.,
 319 (2013), no. 2, 535--555.

\bibitem{NovakBanach} J. Novak \emph{Jucys-Murphy elements and the unitary Weingarten function}. Noncommutative harmonic analysis with applications to probability II, 231--235, Banach Center Publ., 89, Polish Acad. Sci. Inst. Math., Warsaw, 2010.

\bibitem{Novak} J. Novak \emph{Three lectures on free probability theory}. arXiv:1205.2097, to appear in MSRI Publications.


\bibitem{Rider2003} B. Rider \emph{A limit theorem at the edge of a non-Hermitian random matrix ensemble}. 
J. Phys. A: Math. Gen. 36, 3401--3409 (2003).

\bibitem{RiderSinclair} B. Rider, C. Sinclair \emph{Extremal laws for the real Ginibre ensemble}. Ann. Appl. Probab. 24 (2014), no. 4, 1621--1651.

\bibitem{RUDAMSSC} M. Rudelson,  \emph{Lecture notes on non-aymptotic random matrix theory}, notes from the AMS Short Course on Random Matrices, 2013.

\bibitem{RUDICM} M. Rudelson, R. Vershynin \emph{Non-asymptotic theory of random matrices: extreme singular values}, Proceedings of the International Congress of Mathematicians. Volume III, 1576--1602, Hindustan Book Agency, New Delhi, 2010.

 \bibitem{RUD} M. Rudelson, R. Vershynin  \emph{Invertibility of random matrices: unitary and orthogonal perturbations},    J. Amer. Math. Soc.  27 (2014), 293--338.

\bibitem{yinLCL3} I. Yin \emph{The local circular law III: general case}. Probab. Theory Relat. Fields (2014) 160:679--732.
\end{thebibliography}
